\documentstyle[fleqn,12pt]{article}
\begin{document}
\pagenumbering{arabic}
\setcounter{page}{1}
\pagestyle{plain}
\baselineskip=18pt

\begin{center}
{\bf The properties of the quantum supergroup $GL_{p,q}(1|1)$ }
\end{center}
\vspace{0.5cm}

\noindent
Salih Celik \\ 
{\footnotesize Mimar Sinan University, Department of Mathematics, 
80690 Besiktas, Istanbul, TURKEY. }\\
and \\
Sultan A. Celik\\
{\footnotesize Yildiz Technical University, Department of Mathematics, 
Sisli, Istanbul, TURKEY. }

\vspace{0.5cm}
{\bf Abstract}

In this paper properties of the quantum supermatrices in the quantum 
supergroup $GL_{p,q}(1|1)$ are discussed. It is shown that any element of 
$GL_{p,q}(1|1)$ can be expressed as the exponential of a matrix of 
non-commuting elements, like the group $GL_q(1|1)$. An explicit construction 
of this exponential representation is presented. 

\vfill\eject\noindent
{\bf 1. Introduction}

In the past few years, Drinfeld [1] and Faddeev et al [2] constructed a new 
mathematical object, called the quantum group which later has been 
generalized to the quantum supergroup. This topic has been studied by many 
mathematicians and theoretical physicists. 

The simplest supergroup is the group of 2x2 supermatrices with 
two even and two odd matrix elements, i.e. $GL(1|1)$. Even 
matrix elements commute with everything and odd matrix elements 
anticommute among themselves. The deformation of the supergroup 
of 2x2 matrices, i.e. the quantum supergroup $GL_q(1|1)$ can be 
found in Refs. 3-5. A two parameter deformation of $GL(1|1)$ was 
given in Ref. 6 and also Ref. 7. 

It was shown in Ref. 5 that any element of $GL_q(1|1)$ can be 
written as the exponential of a matrix and this exponential 
form was explicitly constructed. This work will be along the lines 
of the work in Ref. 5. In sec. 2 we present a review of $GL_{p,q}(1|1)$ 
together proofs some elementary lemmas which we use in the later 
sections. In sec. 3 we get the matrix elements of $T^n$, the $n$-th 
power of $T$, for $T \in GL_{p,q}(1|1)$. Using these matrix elements 
we prove that $T^n \in GL_{p^n,q^n}(1|1)$ if $T \in GL_{p,q}(1|1)$. 
This result suggests that an element of $GL_{p,q}(1|1)$ can be expressed 
as an exponential of a matrix whose entries obey $(h_1,h_2)$-dependent 
commutation relations (here the parameters $h_1$ and $h_2$ are the logarithm 
of the deformation parameters $q$ and $p$, respectively). To prove this 
in sec. 4 we use the method of the paper of Schwenk et al [5]. We 
derive the explicit form of the $n$-th power of the matrix $M$ which 
is the natural logarithm of the matrix $T \in GL_{p,q}(1|1)$ in sec. 4.  
Thus, we obtain the matrix elements of $T$ in terms of $M$ and vice versa. 
Finally, we state that the usual relation between the superdeterminant and 
the supertrace, which is also satisfy in the supergroup $GL_q(1|1)$, is true 
again in $GL_{p,q}(1|1)$. 
 
\noindent
{\bf 2. Notations and useful formulas } 

In Ref. 3 Manin identifies a quantum supergroup with the 
endomorphisms acting on quantum superplanes. In the matrix 
representation of these endomorphisms, the commutation 
relations of the space coordinates  include the commutation 
relations of the group parameters. 

We state briefly some notations and useful formulas we are going 
to need in this work. The quantum supergroup $GL_{p,q}(1|1)$ consists 
of all matrices in the form 
$$ T = \left ( \matrix{ a & \beta \cr
                         \gamma & d \cr}
\right) \eqno(2.1) $$
where the elements $a, \beta, \gamma$ and $d$ obey the following 
commutation relations 
  $$ a\beta = q \beta a , ~~d\beta = q\beta d, $$
  $$ a\gamma = p\gamma a ,~~ d\gamma = p\gamma d, \eqno(2.2) $$
  $$ \beta \gamma + pq^{-1} \gamma \beta = 0 , ~~ \beta^2 = 0 = \gamma^2, $$
  $$ ad - da = (p - q^{-1})\gamma \beta $$
and $p, q$ non-zero complex numbers, $pq \pm 1 \neq 0$. 

The quantum superdeterminant is defined as 
$$ sD_{p,q}(T) = ad^{-1} - \beta d^{-1} \gamma d^{-1} \eqno(2.3)$$
provided $d$ is invertible. We will also suppose that $a$ is invertible. 
Using (2.2) it is easy to show that $sD_{p,q}(T)$ is central, that is, 
it commutes with $a, a^{-1}, d, d^{-1}, \beta$ and $\gamma$. 

If we take 
$$ \Delta_1 = ad -p^{-1} \beta \gamma \qquad \mbox{and} \qquad 
   \Delta_2 = da -q^{-1} \gamma \beta \eqno(2.4) $$
the super-inverse of $T$ becomes 
$$ T^{-1} = \left ( \matrix{ d\Delta_1^{-1} & -q^{-1} \beta \Delta_2^{-1} \cr
               -p^{-1} \gamma \Delta_1^{-1} & a\Delta_2^{-1} \cr}
\right). \eqno(2.5) $$
Then the superdeterminant is given by 
$$ sD_{p,q}(T) = a^2 \Delta_2^{-1}. \eqno(2.6) $$
Of course, the superdeterminant of $T^{-1}$ may also be defined and it is of 
the form 
$$ sD_{p,q}(T^{-1}) = d^2 \Delta_1^{-1}. \eqno(2.7) $$
On the other hand,
$$ \left \{a^2 \Delta_2^{-1} \right \}^{-1} = d^2 \Delta_1^{-1}. \eqno(2.8) $$
So, 
$$ sD_{p,q}(T^{-1}) = \left \{sD_{p,q}(T) \right \}^{-1}. \eqno(2.9) $$
This relation is generalized in Corollary of sec. 3.

Before passing to the next section, we give three lemmas. 

\noindent
{\bf Lemma 2.1}. For any integer $n$ 
$$ (a - \beta d^{-1} \gamma)^n = a^n - {{q^n - p^{-n}}\over {q - p^{-1}}} 
                                  \beta a^{n-1}d^{-1} \gamma. \eqno(2.10)$$

{\it Proof}. The relation (2.10) can be proved by a induction procedure. 

We note here that if $f$ is any function of $a$ and $d$ of the form 
$f(a,d) = a^n d^m$ where $n$ and $m$ are integer, in the products with $\beta$ 
or $\gamma$, the arguments $a^n$ and $d^m$ of the function $f(a,d)$ behave as 
commuting quantities. So in the products by $\beta$ or $\gamma$, the element 
$a^n$ commutes with $d^m$, i.e., 
$$\beta a^n d^m = \beta d^m a^n.$$

\noindent
{\bf Lemma 2.2}. For any integers $n$ and $m$ 
$$a^n d^m = d^m a^n + (p^n - q^{-n}) { {p^m - q^{-m}} \over {p - q^{-1}}} 
  \gamma a^{n-1} d^{m-1} \beta. \eqno(2.11)$$

{\it Proof}. The proof of this lemma can also be proved by an induction 
procedure. 

\noindent
{\bf Lemma 2.3}. For any integer $n$
$$ \left \{ sD_{p,q}(T) \right \}^n = a^nd^{-n} - 
    p {{p^{-n} - q^n}\over {p - q^{-1}}} a^{n-1} \gamma d^{-n-1} \beta. 
        \eqno(2.12) $$

{\it Proof}. With (2.6), one can write 
$$ \left \{ sD_{p,q}(T) \right \}^n = a^{2n}\Delta_2^{-n}, \eqno(2.13)$$
since $a$ and $\Delta_2$ commute. On the other hand, it may be shown by using 
(2.11) with $n = m$ [or from (2.10)] that 
$$ \Delta_2^n = a^n d^n - p {{p^n - q^{-n}}\over {p - q^{-1}}} 
                  a^{n-1} \gamma d^{n-1} \beta.$$
Hence, by replacing $n$ by $-n$ into the above equation one gets 
$$ a^{2n}\Delta_2^{-n} = a^n d^{-n} - p {{p^{-n} - q^n}\over {p - q^{-1}}} 
   a^{n-1)} \gamma d^{-n-1} \beta = \left \{sD_{p,q}(T) \right \}^n, $$
as required.

These results will be used in the following sections. 

\noindent
{\bf 3. The properties of $T^n$ }

To show that $T^n \in GL_{p^n,q^n}(1|1)$ for $T \in GL_{p,q}(1|1)$ 
we will explicitly obtain the matrix elements of $T^n$, the $n$-power 
of $T$, for a matrix $T \in GL_{p,q}(1|1)$. The matrix elements of 
$T^n$ in a more compact form also appear in the paper of Schwenk et al [5]. 
First we define the following functions. Let 
$$F_n(a,q^{-1}d) \beta \gamma = 
     \sum_{k=0}^{n-2} <n - k - 1>_{pq} a^{n-k-2}(q^{-1}d)^k \beta \gamma, 
      \eqno(3.1)$$
$$G_n(a,q^{-1}d) \beta = \sum_{k=0}^{n-1} a^{n-k-1}(q^{-1}d)^k \beta, 
                           \eqno(3.2)$$
where
$$ <N>_{pq} = {{1 - (pq)^{-N}}\over {1 - (pq)^{-1}}}. \eqno(3.3) $$

\noindent
{\bf Lemma 3.1.} If $T \in GL_{p,q}(1|1)$ then the matrix $T^n$, the $n$-th 
power of $T$, has the form 
$$ T^n = \left ( \matrix{ A_n & B_n \cr
                          C_n & D_n \cr}
\right) \eqno(3.4) $$
where 
$$ A_n = a^n + F_n(a,q^{-1}d) \beta \gamma, ~~ B_n = G_n(a,q^{-1}d) \beta, $$
$$ D_n = d^n + F_n(d,p^{-1}a) \gamma \beta, ~~C_n = G_n(d,p^{-1}a) \gamma. 
\eqno(3.5)$$

{\it Proof.} It can be done by induction on $n$ using the fact that 
$T^{n+1} = T^n T$. 

\noindent
{\bf Lemma 3.2.} If $T \in GL_{p,q}(1|1)$ then 
$T^n \in GL_{p^n,q^n}(1\vert 1)$. That is, the matrix elements of $T^n$ obey 
the following commutation relations 
$$ A_n B_n = q^n B_n A_n, \qquad D_n B_n = q^n B_n D_n,$$
$$ A_n C_n = p^n C_n A_n, \qquad D_n C_n = p^n C_n D_n, \eqno(3.6)$$
$$ B_n^2 = 0 = C_n^2, \qquad q^n B_n C_n + p^n C_n B_n = 0,$$
and 
$$ [A_n, D_n] = (p^n - q^{-n}) C_n B_n, \eqno(3.7)$$
where 
$$ [u,v] = uv - vu. $$

{\it Proof.} It is not difficult to check that the relations (3.6) are 
satisfied. The reader can easily prove this by using the relations (2.2) and 
the equations (3.5). But proof of the relation (3.7) requires some operations. 
The proof of (3.7) can be found in the Appendix. 

The property in the Lemma 3.2 gives an opportunity for the matrix 
$T \in GL_{p,q}(1|1)$ to be represented as an exponential of a matrix. 
This will be done in the last section. 

Now we want to calculate the superdeterminant of $T^n$. For this, we 
write the matrix $T^n$ in the form 
$$ T^n = \left ( \matrix{ A_n & 0 \cr
                          C_n & D_n - C_nA_n^{-1}B_n \cr}
\right) 
         \left ( \matrix{ 1 & A_n^{-1}B_n \cr
                          0 & 1 \cr}
\right)
\eqno(3.8) $$
using the Crout decomposition. Then the superdeterminant of $T^n$ becomes 
$$sD_{p,q}(T^n) = A_n \left (D_n - C_nA_n^{-1}B_n \right )^{-1} 
   = \left (A_n - B_nD_n^{-1}C_n \right )D_n^{-1}.\eqno(3.9)$$
After some calculations one gets 
$$ sD_{p,q}(T^n) = a^nd^{-n} - p {{p^{-n} - q^n}\over {p - q^{-1}}} 
                   a^{n-1} \gamma d^{-n-1} \beta \eqno(3.10) $$ 
which is the same with (2.12). Thus we have: 

\noindent
{\bf Corollary 3.3.}. For any integer $n$ 
$$ \left \{ sD_{p,q}(T) \right \}^n  =  sD_{p,q}(T^n). \eqno(3.11) $$ 

\noindent
{\bf 4. The exponential parametrization of $GL_{p,q}(1|1)$ }

The implication $T^n \in GL_{q^n}(1|1)$ for $T \in GL_q(1|1)$ 
suggests that $T$ can be represented by exponentiating a matrix. This was 
shown by Schwenk et al [5]. However, in sec. 3 it is show that 
$T^n \in GL_{p^n,q^n}(1|1)$ if $T \in GL_{p,q}(1|1)$ [see, Lemma 3.2]. Thus 
the matrix $T \in GL_{p,q}(1|1)$ can be written as the exponential of a matrix 
with non-commuting entries. Let 
$$ q = e^{h_1} \qquad \mbox{and} \qquad p = e^{h_2}. \eqno(4.1)$$
Suppose that
$$ T = e^{hM}, \qquad h = {{h_1 + h_2} \over 2} \eqno(4.2)$$
where
$$ M = \left ( \matrix{ x & \mu \cr
                        \nu & y \cr}
\right). \eqno(4.3) $$

To find the commutation relations of the matrix elements of $M$ we write 
the exponent as 
$$ M = {1\over h}\ln T. \eqno(4.4)$$
The logarithm of the matrix $T$ is defined by 
$$ \ln T = \sum_{n=1}^{\infty} {(-1)^{n+1}\over n} (T - I)^n \eqno(4.5) $$
as a series expansion. Now let 
$$ (T - I)^n = \left ( \matrix{ \tilde{A}_n & \tilde{B}_n \cr
                              \tilde{C}_n & \tilde{D}_n \cr}
\right). \eqno(4.6) $$
Then some calculations show that 
$$ \tilde{A}_n = (a - 1)^n + \sum_{j=0}^{n-2}\sum_{k=0}^{n-j-2} (a - 1)^k 
                 (p^{-1}q^{-1}a - 1)^{n-k-j-2} (q^{-1}d - 1)^j \beta \gamma, $$
$$ \tilde{D}_n = (d - 1)^n + \sum_{j=0}^{n-2}\sum_{k=0}^{n-j-2} (d - 1)^k 
                 (p^{-1}q^{-1}d - 1)^{n-k-j-2} (p^{-1}a - 1)^j \gamma \beta, 
                  \eqno(4.7)$$
$$ \tilde{B}_n = \sum_{j=0}^{n-1} (a - 1)^{n-j-1} (q^{-1}d - 1)^j \beta, \qquad
   \tilde{C}_n = \sum_{j=0}^{n-1} (d - 1)^{n-j-1} (p^{-1}a - 1)^j \gamma. $$
We want to obtain the matrix elements of $M$ in an explicit form. For the sake 
of simplicity we define 
\begin{eqnarray*}
f_q(a,d,p) \beta \gamma 
     & = & {{q^2}\over {q - p^{-1}}} \left ( {{\ln a}\over {a(qa - d)}} - 
            {{\ln (p^{-1}q^{-1}a)}\over {a(p^{-1}a - d)}} 
            \right )\beta \gamma  \\ 
     &   & + q^2 {{\ln (q^{-1}d)}\over {(p^{-1}a - d)(qa - d)}} 
            \beta \gamma \hspace*{4.3cm}{(4.8)}
\end{eqnarray*}
$$ g(a,q^{-1}d) \beta = {{\ln a - \ln (q^{-1}d)}\over {a - q^{-1}d}} \beta, 
  \eqno(4.9)$$
where the logarithms of $a$ and $d$ exist and the denominators are non-zero. 
Then we have 

\noindent
{\bf Lemma 4.1.} If $T \in GL_{p,q}(1\vert 1)$ then the expressions for the 
matrix elements of $M$ in terms of $T$ are as folows: 
$$ x = {1\over h} \left \{ \ln a + f_q(a,d,p) \beta \gamma \right \}, \qquad 
   \mu = {1\over h} g(a,q^{-1}d) \beta, $$
$$ y = {1\over h} \left \{ \ln d + f_p(d,a,q)\gamma \beta \right \}, \qquad 
   \nu = {1\over h} g(d,p^{-1}a) \gamma. \eqno(4.10) $$

{\it Proof.} Use the equation (4.4) with (4.5). 

Note that the matrix elements $a$ and $d$ are behave as commuting 
quantities when they are in a product case by $\beta$ or $\gamma$. 
Thus it is not necessary to order the arguments in the equations (4.8) and 
(4.9). 

\noindent
{\bf Proposition 4.2.} If $T \in GL_{p,q}(1\vert 1)$ then the matrix elements 
of $M$ obey the following commutation relations 
$$ [x,\mu] = {{2h_1}\over {h_1 + h_2}} \mu, \qquad 
   [y,\mu] = {{2h_1}\over {h_1 + h_2}} \mu, \qquad \mu^2 = 0,$$
$$ [x,\nu] = {{2h_2}\over {h_1 + h_2}} \nu, \qquad 
   [y,\nu] = {{2h_2}\over {h_1 + h_2}} \nu, \qquad \nu^2 = 0, \eqno(4.11) $$
$$  xy - yx = 0, \qquad \mu \nu + \nu \mu = 0. $$

{\it Proof.} It is easy to see that the relations (4.11) [if we except that 
the relation $xy - yx = 0$] is satisfied. Let us prove the last relation in 
(4.11), here. Let 
$$X = [\ln a, \ln d],$$
$$Y = [\ln a, f_p(d,a,q)\gamma\beta],$$
$$Z = [\ln d, f_q(a,d,p)\beta\gamma].$$
After some calculations one gets
$$Y - Z = {{4 h^2}\over {1 - pq}} \gamma a^{-1} \beta d^{-1}, \eqno(4.12a)$$
and 
$$X = {{\ln^2 (pq)}\over {pq - 1}} \gamma a^{-1} \beta d^{-1}. \eqno(4.12b)$$
[The proof of (4.12b) is rather lengthy but straightforward]. Thus we have 
$$X + Y - Z = 0.$$
Here we used the relations 
$$[\ln a,\beta] = h_1 \beta = [\ln d,\beta], $$
$$ [\ln a,\gamma] = h_2 \gamma = [\ln d,\gamma]. \eqno(4.13)$$

If $T \in GL_q(1|1)$ and $T = e^{\theta M}$ where the matrix $M$ is given by 
(4.3) and $q = e^\theta$, we know that $x - y$ is the central element, which 
is known as the supertrace of the matrix $M$. This case is true in the 
supergroup $GL_{p,q}(1|1)$ too, i.e. the supertrace of $M$, 
$str M = x - y$, is the central element in the algebra (4.11).

Now we will obtain the matrix elements of $T$ in terms of $M$. First 
we derive the explicit form of $M^n$, the $n$-th power of $M$. It will 
simplify the elements of $M^n$ to define a transformation $\tau$ by 
$$ \tau: \tau(x,y,\mu,\nu,h_1,h_2) \mapsto (y,x,\nu,\mu,h_2,h_1). \eqno(4.14)$$
Then the relation (4.11) are preserved by $\tau$. For example, 
$$ [x,\mu]^\tau = ({{h_1}\over h} \mu)^\tau ~\Longrightarrow~ 
   [y,\nu] = {{h_2}\over h} \nu. $$
Let 
$$\phi = {{h_1}\over h} ~~\mbox{and}~~ \varphi = {{h_2}\over h}. \eqno(4.15)$$

The following lemma can be proved by mathematical induction. We denote the 
algebra (4.11) by ${\cal M}_{h_1,h_2}$. 

\noindent
{\bf Lemma 4.3.} If $M \in {\cal M}_{h_1,h_2}$ then the matrix $M^n$, has the 
form 
$$ M^n = \left ( \matrix{ x^n - \mu \nu F_n & \mu G_n \cr
                          \nu G_n^\tau      & y^n - \nu \mu F_n^\tau \cr}
\right) \eqno(4.16) $$
where
\begin{eqnarray*}
F_n & = & F_n(x,y,\phi,\varphi) \\
    & = & {x^n \over {2(x - y - \varphi)}} - 
          {(x + \phi + \varphi)^n \over {2(x - y + \phi)}} - 
          {(y + \varphi)^n \over {(x - y + \phi)(x - y - \varphi)}} 
\hspace*{0.3cm}{(4.17)}
\end{eqnarray*}
\begin{eqnarray*}
G_n & = & G_n(x,y,\phi) \\
    & = & {{(x + \phi)^n - y^n} \over {x - y + \phi}}. \hspace*{7.7cm}{(4.18)}
\end{eqnarray*}

Now we easily obtain the  expressions for the matrix elements of $T$ in terms 
of $M$ using the equation (4.2): 

\noindent
{\bf Lemma 4.4.} If $T = e^M$ then one has 
$$ a = e^{hx} - {{\mu \nu}\over {(x - y + \phi)(x - y - \varphi)}}\left \{
       \left( {{\phi + pq\varphi}\over 2} - {{pq - 1}\over 2}(x - y) \right)
       e^{hx} - pe^{hy} \right \},$$
$$ d = e^{hy} - {{\nu \mu}\over {(x - y + \phi)(x - y - \varphi)}}\left \{
       \left ( {{\varphi + pq\phi}\over 2} + {{pq - 1}\over 2}(x - y) \right )
       e^{hy} - qe^{hx} \right \}, $$
$$ \beta = {{\mu}\over {x - y + \phi}}(qe^{hx} - e^{hy}), ~~
 \gamma = {{\nu}\over {\varphi - (x - y)}}(pe^{hy} - e^{hx}). \eqno(4.19)$$

\noindent
{\bf Proposition 4.5.} If $M \in {\cal M}_{h_1,h_2}$ and $T = e^M$ then 
$T \in GL_{p,q}(1\vert 1)$. 

{\it Proof.} To prove that the matrix $T$ is in $GL_{p,q}(1\vert 1)$ the 
reader can be verified the relations (2.2). 

In Ref. 5 it has been shown that the usual relation between the 
superdeterminant and supertrace is valid, i.e.
$$ sD_q(T) = e^{\theta strM}, \qquad q = e^{\theta}. $$ 
Finally with direct calculation we can show that 
$$ sD_{p,q}(T) = e^{h(x - y)} = e^{h strM}, \eqno(4.20)$$
where $h = {1\over 2} \ln (pq)$. 

{\it Remarks}. In the equations (4.11) if we take $h_1 = h_2$ we obtain 
the algebra in Ref. 5 (the equ.s (5.7)) where $q$ replaces $q^{-1}$ for 
$p = q$. In this case the relations (4.19) identify the equations (5.9) in 
Ref. 5. Thus our work can be considered as a generalization of Ref. 5. 

\noindent
{\bf Appendix: the proof of eq. (3.7)}

Now we will show that
$$ [A_n,D_n] = (p^n - q^{-n})C_nB_n. \eqno(A1) $$
For this, we will use the fact that
$$ T^{k + 1} = T^kT = TT^k. \eqno(A2) $$
It is proved by induction on $n$. 

(1) For $n = 1$, the equality (A1) identifies with the last relation in (2.2).

(2) Assume that the equation is true for $n = k$.

(3) With (A2) we write that
$$ A_{k + 1} = A_1A_k + B_1C_k, \qquad C_{k + 1} = C_1A_k + D_1C_k, $$
$$ B_{k + 1} = A_1B_k + B_1D_k, \qquad D_{k + 1} = D_1D_k + C_1B_k. $$
Now some calculations show that
\begin{eqnarray*}
[A_{k + 1},D_{k + 1}] & = & \left ( (pq)^{k+1} - 1 \right ) 
                            \left ( C_1B_kA_1A_k - (pq)^{-k-1} D_1D_kB_1C_k 
                             \right ) + K \\
                      & = & \left (p^{k+1} - q^{-k-1} \right )C_{k + 1}B_{k+1} 
                             + K - L, 
\end{eqnarray*}
where 
\begin{eqnarray*}
 K & = & \left (p^{2k+1}q^k - q^{-k-1} \right )C_1A_kB_1D_k + 
         \left ( p^{k+1} - pq^{-k} \right )C_kA_1B_kD_1 + \\
   &   & p^{k+1} \left ( C_kA_1A_kB_1 - C_1A_kA_1B_k \right ) + \\
   &   & p^{k+1} \left ( D_kC_1B_kD_1 - D_1C_kB_1D_k \right ) 
\end{eqnarray*}
and 
\begin{eqnarray*}
 L & = & \left (p^{k+2}q - pq^{-k} \right ) C_kA_1B_kD_1 + 
            \left (p^{k+1} - q^{-k-1} \right ) C_1A_kB_1D_k. 
\end{eqnarray*}
Thus it must be 
$$ K - L = 0. $$
In fact, 
\begin{eqnarray*}
K - L & = & p^{k+1} \left (p^kq^k - 1 \right ) C_1A_kB_1D_k + 
            p^{k+1} \left (1 - pq \right ) C_kA_1B_kD_1 + \\
      &   & p^{k+1} \left ( C_kA_1A_kB_1 - C_1A_kA_1B_k \right ) + \\
      &   & p^{k+1} \left ( D_kC_1B_kD_1 - D_1C_kB_1D_k \right ) \\
      & = & 0.
\end{eqnarray*}

\end{document}